\documentclass[12pt,reqno]{article}

\usepackage[usenames]{color}
\usepackage{amssymb}
\usepackage{amsmath}
\usepackage{amsthm}
\usepackage{amsfonts}
\usepackage{amscd}
\usepackage{graphicx}

\usepackage[colorlinks=true,
linkcolor=webgreen,
filecolor=webbrown,
citecolor=webgreen]{hyperref}

\definecolor{webgreen}{rgb}{0,.5,0}
\definecolor{webbrown}{rgb}{.6,0,0}

\usepackage{color}
\usepackage{fullpage}
\usepackage{float}

\usepackage{graphics}
\usepackage{latexsym}
\usepackage{epsf}
\usepackage{breakurl}

\newcommand{\ubar}[1]{\text{\b{$#1$}}}

\newcommand{\seqnum}[1]{\href{https://oeis.org/#1}{\rm \underline{#1}}}
\def\modd#1 #2{#1\ \mbox{\rm (mod}\ #2\mbox{\rm )}}

\DeclareMathOperator{\lcm}{lcm}

\begin{document}
	
	\begin{center}
		\epsfxsize=4in
		\leavevmode
		
	\end{center}
	
	\theoremstyle{plain}
	\newtheorem{theorem}{Theorem}
	\newtheorem{corollary}[theorem]{Corollary}
	\newtheorem{lemma}[theorem]{Lemma}
	\newtheorem{proposition}[theorem]{Proposition}
	
	\theoremstyle{definition}
	\newtheorem{definition}{Definition}
	\newtheorem{example}{Example}
	\newtheorem{conjecture}[theorem]{Conjecture}
	
	\theoremstyle{remark}
	\newtheorem{remark}[theorem]{Remark}

	\begin{center}
		\vskip 1cm{\LARGE\bf Rigged Horse Numbers and their Modular Periodicity\\
			\vskip .1in
		}
		\vskip 1cm
		\large
		Benjamin Schreyer\\
		
		Departments of Computer Science and Physics\\
		University of Maryland \\
		College Park, Maryland 20742\\
		USA\\
		Plasma Physics Division\\
		U.S. Naval Research Laboratory\\
		Washington, D.C. 20375\\
		USA\\
		
		\href{mailto: benontheplanet@gmail.com}{\tt benontheplanet@gmail.com} \\
		
	\end{center}
	
	\vskip .2 in
	
	\begin{abstract}
		The {\em Fubini numbers} count the permutations of horse racing where ties are possible. The closely related {\em $r$-horse numbers} count the finishes of a horse race where some subset of $r$ horses agree to finish the race in a specific relative strong ordering. We express the $r$-Fubini numbers as a sum of $r$ index-shifted sequences of Fubini numbers weighted with the {\em signed Stirling numbers of the first kind}. We use a novel shift operator counting. Further, we demonstrate the eventual modular periodicity of {\em $r$-Fubini numbers}. Their maximum period is determined to be the {\em Carmichael function} of the modulus. The maximum period occurs in the case of an odd modulus for Fubini numbers.
	\end{abstract}
	
	\section{Introduction}
	\subsection{Orderings weak and strong}
		The notation $\mathbb{N}$ denotes the set of positive integers. The symbol $\mathbb{N}_{0} = \mathbb{N} \cup \{0\}$ represents the set of nonnegative integers.
	\begin{definition}
		{\em Fubini numbers}, denoted $F(n)$, count weak orderings of $n$ elements, with $n \in \mathbb{N}_{0}$.
	\end{definition}
	In the case of a horse race, the ordering is weak, an equivalence determines a tie, and strict less-than or greater-than relations determine a clear succession of horses. The horse numbers or ordered Bell numbers are other common names for the number of weak orderings. When ties are impossible, these orderings become regular (strong) permutations. Velleman and Call \cite{cc:velleman} gave another intuitive example for Fubini number counting: the number of combinations for a combination lock in which all buttons are used exactly once, and multiple buttons can be pressed simultaneously.
	\subsection{Rigged weak orderings}\label{ss:riggedweakorderings}
	\begin{definition}
		The constraint of {\em relative strong ordering} applies to a subset $X'$ of a weakly ordered set $X$. The constraint demands that relations between elements of $X'$ are only strict less-than or greater-than relations. An element $x' \in X'$ can be set equivalent to another $x$ only if $x \notin X'$.
	\end{definition}
	\begin{definition}
		The {\em $r$-Fubini numbers}, $F_{r}(n)$, count weak orderings of a set with cardinality $n$ such that $r$ elements are distinguished and constrained to follow relative strong ordering.
	\end{definition}
	\begin{definition}
		The {\em $r$-horse numbers}, $H_{r}(n)$, count weak orderings such that $r$ elements of a set of cardinality $n$ are distinguished and constrained to follow a specified strong permutation relative to each other. 
	\end{definition}
	Every permutation counted by $H_{r}(n)$ has $r!$ rearrangements of the specified permutation $x_{1} < x_{2} < \cdots < x_{r}$. The rearrangements cover all cases where the $r$ elements are mutually inequivalent. We can write this as the equation
	\begin{align}
		r! H_{r}(n) = F_{r}(n).\label{eqn:HFrelation}
	\end{align}
	The $r$-horse numbers count a nontrivial restriction on weakly ordered permutations. A restriction where two or more elements must be tied is the same as reducing the effective number of elements in the counting. The equivalent elements act as a single unit in every permutation. 
	\subsection{The Stirling numbers of the first and second kind}\label{sec:stirling}
	\begin{definition}
		{\em The signed Stirling numbers of the first kind}, $s(n,k)$, count partitions of $n$ elements into $k$ cycles. The sign gives the parity of $n - k$. For combinatorics, we have $n,k \in \mathbb{N}_{0}$. The two-index sequence $s(n,k)$ can be arranged into a matrix. We index rows with $n$ and columns with $k$. Let $\hat{s}$ denote the infinite matrix of $s(n,k)$.
	\end{definition}
	\begin{definition}
		{\em The Stirling numbers of the second kind}, $S(n,k)$, count ways to partition a set of $n$ elements into $k$ subsets. The matrix of $S(n,k)$, labeled $\hat{S}$, uses the same indexing as $\hat{s}$.
	\end{definition}
	We introduce the Stirling numbers with the addition of three useful properties. The notation $I$ is the identity operation. \textit{Advanced Combinatorics} \cite[Eqn.\ 6f,\ p.\ 144]{cc:matrix} gives the first proposition.
	\begin{proposition}
		The matrices $\hat{s}$ and $\hat{S}$ are inverses of each other, so 
		\begin{align}
			\hat{s} \hat{S} = \hat{S} \hat{s} = I \label{eqn:stirlinginvs}.
		\end{align}
		Both $\hat{s}$ and $\hat{S}$ are lower triangular.
	\end{proposition}
	{\em Concrete Mathematics} \cite[Eqn.\ 6.13,\ p.\ 263]{cc:cm} provides the second important property.
	\begin{proposition}\label{prop:fallfac}
		The Stirling numbers of the first kind give the coefficients for powers of the argument of the falling factorial, so
		\begin{align}
			x (x - 1) \cdots (x - n + 1) = \sum_{k = 0}^{n} s(n,k) x^{k}\label{eqn:fallfac}.
		\end{align}
	\end{proposition}
	The falling factorial of $x$ with $n$ multiplicative terms, $x (x - 1) \cdots (x - n + 1)$, is written more succinctly as $(x)_{\ubar{n}}$.
	We give a definition needed to discuss eventual modular periodicity.
	\begin{definition}
		A sequence $f(n)$ is {\em eventually periodic modulo $K$}, with $K \in \mathbb{N}$, if there exists a $q \in \mathbb{N}$ such that
		for large enough $a \in \mathbb{N}_{0}$:
		\begin{align}
			f(a)~{} \equiv f(a + q) \pmod{K}.
		\end{align}
	\end{definition}
	Finally, $S(n,k)$, with fixed $k$, is eventually modular periodic in $n$ under any modulus. The periodicity property can be shown using the following formula from the paper \textit{Stirling matrix via Pascal matrix} \cite[p.\ 55]{cc:relation}.
	\begin{lemma}\label{lem:sskexplicit}
		The Stirling numbers of the second kind have the explicit form
		\begin{align}
			S(n,k) = \frac{1}{k!} \sum_{t = 0}^{k}(-1)^{k-t}\binom{k}{t} t^{n}\label{eqn:stirlingsecond}.
		\end{align}
	\end{lemma}
	The periodicity of modular exponentiation determines eventual modular periodicity for $S(n,k)$ for fixed $k$.
	\subsection{The Carmichael function}
	\begin{definition}
		The function $\lambda(K)$ is the {\em Carmichael function}. The Carmichael function gives the least common multiple of all periods of integer exponentiation modulo $K \in \mathbb{N}$. It is often useful in the context of the multiplicative group of integers modulo $K$.
	\end{definition}
	The Carmichael function has the following recurrence \cite{cc:carmichael} using {\em Euler's totient function} $\varphi(n)$:
	\begin{equation}
		\lambda(K) = \begin{cases}
			\varphi(K), & \text{if $K \in \{1,2,4\}$ or $n$ is an odd prime power;}\\
			\frac{1}{2}\varphi(K), & \text{if $K = 2^{r}, r \geq 3$;}\\
			\lcm (\lambda(p_{1}), \lambda(p_{2}), \ldots), & \text{if $K = p_{1}^{a_{1}}p_{2}^{a_{2}} \cdots p_{R}^{a_{R}}$.}
		\end{cases}\label{eqn:carmichaelrecur}
	\end{equation}
 We state two properties of modular exponentiation \cite[p.\ 190]{cc:jsbc}.
	\begin{proposition}
		For all $a$, $a \in \{0, 1, \ldots, K - 1\}$, we have
		\begin{align}
			a^{R} \equiv a^{\lambda(K) + R} \pmod {K}.
		\end{align}
		Here $R = \max(R_{1}, R_{2}, \ldots, {R_{N}})$ given $K = p_{1}^{R_{1}} p_{2}^{R_{2}}\cdots p_{N}^{R_{N}}$.
	\end{proposition}
	For integers coprime to $K$ a stronger statement holds.
	\begin{proposition}
		For all $b$ coprime to $K$
		
		\begin{align}
			b^{\lambda(K)} \equiv 1 \pmod {K}.
		\end{align}
	\end{proposition}
	\subsection{Operations preserving eventual periodicity}\label{subsec:evperioprop}
	Scaling by an integer, addition with another modular periodic sequence, and index shifting preserve the eventual modular periodicity of a sequence. Upper bounds for the period do not change under scaling and shifting. The least common multiple of the sequences' periods must be considered to include the addition of sequences.
	\subsection{Shift operators}
	We use shift operators to formally show that $r$-Fubini numbers are expressible using the signed Stirling numbers of the first kind.
	\begin{definition}
		Operators $E$ and $E^{-1}$ are the left and right shift operators, respectively. We abbreviate repeated shift operations as $E^{m}$ with $m\in \mathbb{Z}$. A zero shift $E^{0}$ is also the identity operation $I$.
	\end{definition}
	Computation of the $r$-Fubini numbers uses the left and right shift operators on the sequence $F(0), F(1), \ldots, F(n + r)$. We only consider one-sided sequences. Shift operators are linear when applied to a sequence.  Importantly, sequences also distribute over addition of shift operators, which means $(AE^{a} + BE^{b})C(n) = AC(n + a) + BC(n+b)$. If a sequence of numbers occupies the entries of a vector, a shift operator has a matrix representation. Zero is placed in the first index of the vector sequence when $E^{-1}$ is applied. We introduce vector notation with the upper arrow, indexing by subscript, and provide explicit definitions:
	\begin{align}
		(E\vec{C})_{n} &= \vec{C}_{n + 1}\\
		(E^{-1}\vec{C})_{n} &= \begin{cases}
			0, & \text{if $n = 0$;}\\
			\vec{C}_{n - 1}, & \text{if $n \neq 0.$}
		\end{cases}
	\end{align}
	\subsection{Related works}
	In this text, the $r$-horse and $r$-Fubini counted permutations are indexed by denoting the total number of ordered elements $n$ and the size of the subset that must follow a fixed or arbitrary strong ordering $r$. Other authors \cite{cc:racz} consider a total of $n + r$ elements. R\'acz \cite{cc:racz} studied related orderings and $r$-Fubini numbers, referencing an expression in terms of the $r$-Stirling numbers of the second kind and factorials for $F_{r}(n)$. Broder \cite{cc:broder} characterized the $r$-Stirling numbers and explored related orthogonality relations. Mez\H o gave a proof of periodicity for $F_{r}(n) \pmod{10}$ in a paper on last digit periodicity \cite{cc:mezo}, extending a proof by Gross \cite{cc:gross}. Further, Asgari and Jahangiri \cite{cc:asgari} showed the eventual periodicity of the $r$-Fubini numbers modulo an arbitrary natural number. Asgari and Jahangiri also gave a formula for the period.
	\subsection{Contributions}
	We present a new counting proof. The proof shows that a linear composition of shifted Fubini number sequences is the $r$-Fubini numbers for fixed $r$. The resulting formula is identical in structure to a formula for the {\em $r$-Bell numbers}, which count $r$-partitions, derived by Nyul and Keresk\'enyi-Balogh \cite[Thm.\ 4.2]{cc:nyul}. We further formulate proofs for the eventual modular periodicity of Fubini and $r$-Fubini numbers, which give an upper bound for their modular eventual period. The upper bound is the Carmichael function $\lambda(K)$. The natural number $K$ is the modulus. When $K$ is odd, we show that $\lambda(K)$ is the period. The results apply to the combinatorial problem of $r$-horse numbers under a division by $r!$.
	\section{General rigged orderings of $r \leq n$ elements}	
	 Let $G(n)$ be the number of orderings of a set that contains $n$ elements. The orderings counted observe the restriction of no ties (relative strong ordering) on some subset of $m$ of the $n$ elements. We first establish a lemma required to prove a weighted sum expression for $F_{r}(n)$.
	\begin{lemma}\label{lm:counting}
		Consider a set of $n$ elements with $G(n)$ allowed orderings. Including a new element that is inequivalent with some subset of $m$ elements following relative strong ordering gives $G(n + 1) - mG(n) = (E - mI)G(n)$ orderings.		
		\begin{proof}
			Consider adding the new element with no restriction. The new number of orderings is $G(n + 1)$, since the new element has no constraint upon it. The number of elements is simply increased.
			Restricting the counting to permutations with $x'$ inequivalent to any of the $m$ elements requires the multiplication principle.
			
			When we introduce $x'$ it can be set equivalent to one of $m$ elements to form a disallowed permutation. The choices of the permutation of the original set, counted by $G(n)$, and which equivalence is made on $x'$, counted by $m$, are independent. Independence is clear since all pairs of elements in the relatively strongly ordered subset are inequivalent. Therefore, the multiplication principle determines that $mG(n)$ new disallowed permutations exist. The lemma follows from the complement principle.
		\end{proof}
	\end{lemma}
	The counting above is repeatable for an increasingly large subset that follows a relative strong ordering. The application of each counting step yields a new counting that is compatible with the lemma.
	\subsection{The $r$-Fubini counting with shift operators}
	\begin{theorem}\label{thm:fubinir}
		For indices $n$ and $r$, with $n,r \in \mathbb{N}_{0}$ and $r \leq n$:
		\begin{align}
			H_{r}(n) = \frac{1}{r!} \sum_{j=0}^{r} s(r,r-j) F(n - j).
		\end{align}
		\begin{proof}
			The proof first counts the case where the subset $\{x_{1}, x_{2}, x_{3}, \ldots, x_{r}\}$ follows relative strong ordering, which gives the $r$-Fubini numbers. We divide by $r!$ to specify a permutation of the $r$ elements to get the $r$-horse numbers. To begin counting, we first remove the counting of the $r$ elements to be re-added following relative strong ordering. This leaves weak permutations of $n - r$ elements, which is $E^{- r}F(n)$. Elements of the distinguished subset rejoin the ordered set, with no ties within the subset $\{x_{1}, x_{2}, x_{3}, \ldots, x_{r}\}$. The result is subtractions of ascending integer multiples of the identity operation from $E$ by Lemma \ref{lm:counting}. We make the following $r$ steps, adding back each element of $\{x_{1}, x_{2}, x_{3}, \ldots, x_{r}\}$.
			\begin{itemize}
				\item{Add $x_{1}$ inequivalent to any $x \in \emptyset$. The count is $EE^{-r}F(n)$.}
				\item{Add $x_{2}$ inequivalent to any $x \in \{x_{1}\}$. The count is $(E - 1I)EE^{- r}F(n)$. }
				\item{Add $x_{3}$ inequivalent to any $x \in \{x_{1}, x_{2}\}$. The count is $(E - 2I)(E - 1I)EE^{-r}F(n)$.}
				\item{$\cdots$}
				\item{Add $x_{r}$ inequivalent to any $x \in \{x_{1}, x_{2}, \ldots ,x_{r - 1}\}$. The count is $(E - (r - 1)I)\cdots(E - 2I)(E - 1I)EE^{-r}F(n)$.}
			\end{itemize}
			All elements are included with their desired ordering with those from $\{x_{1}, x_{2}, \ldots, x_{r}\}$ inequivalent. All counted orderings have $\{x_{1}, x_{2}, \ldots, x_{r}\}$ relatively strongly ordered. The falling factorial appears with argument $E$ and $r$ terms, so
			\begin{align}
				F_{r}(n) = (E)_{\ubar{r}} E^{-r} F(n).\label{eqn:fallfacoperator}
			\end{align}
			The falling factorial expands as a sum according to Proposition \ref{prop:fallfac}, giving the equation
			\begin{align}
				F_{r}(n) = \left( \sum_{j = 0}^{r} s(r,j)E^{j} \right) E^{-r} F(n).
			\end{align}
			The proposition applies because the multiplication of operators acts the same as multiplying polynomial variables. The effect of the shift operators is now trivial upon $F(n)$, so $F_{r}(n)$ can be written without operators as
			\begin{align}
				F_{r}(n) = 	\sum_{j = 0}^{r} s(r,j) F(n - r + j).
			\end{align}
			One can re-index the sequences to
			\begin{align}
				F_{r}(n) = 	\sum_{j = 0}^{r} s(r,r - j) F(n - j) \label{eqn:rfubfinal}.
			\end{align}
			Given $F_{r}(n)$, dividing by $r!$ gives $H_{r}(n)$ the form
			\begin{align}
				H_{r} (n) = \frac{1}{r!}\sum_{j = 0}^{r} s(r,r - j) F(n - j) \nonumber.
			\end{align}
		\end{proof}
	\end{theorem}
	An interesting notation can be defined. We place a shift operator in the binomial coefficient. Combining this notation, Equation \ref{eqn:HFrelation}, and Equation \ref{eqn:fallfacoperator} gives
	\begin{align}
		H_{r} (n) &= \binom{E}{r} F(n - r). \label{eqn:binomform}
	\end{align}
	It is interesting to consider the implications of applying such an operator multiple times to a sequence.
	\subsection{A useful alternating recurrence}
	\begin{corollary}\label{cor:auseful}
		For $n\in \mathbb{N}_{0}$, Fubini numbers have the recurrence
		\begin{align}
			F(n) = n! - \sum_{j = 1}^{n} s(n,n - j) F(n - j).
		\end{align}
		\begin{proof}
			If the weakly ordered set is $\{x_{1}, x_{2},\ldots, x_{n}\}$ is constrained such that $\{x_{1}, x_{2},\ldots, x_{n}\}$ are relatively strongly ordered, the number of arrangements is $F_{n}(n)=n!$. It then follows from Theorem \ref{thm:fubinir} that
			\begin{align}
				n! = \sum_{j = 0}^{n} s(n,n - j) F(n - j) \label{eqn:matrixpredec}.
			\end{align}
			We rearrange after the substitution $s(0,0) = 1$ to the result
			\begin{align}
				F(n) = n! - \sum_{j = 1}^{n} s(n,n - j) F(n - j) \nonumber.
			\end{align}
		\end{proof}
	\end{corollary}
	\section{Linear transformation between strong and weak orderings}
	\begin{corollary}\label{cor:lint}
		Let $\vec{f}$ and $\vec{F}$ be infinite vectors with entries $n!$ and $F(n)$, respectively. Vectors $\vec{f}$ and $\vec{F}$ obey the following relations involving Stirling matrices $\hat{s}$ and $\hat{S}$:
		\begin{align}
			\hat{S}\vec{f} &= \vec{F}\label{eqn:strongstirling} \\
			\hat{s}\vec{F} &= \vec{f}\label{eqn:weakstirling}.
		\end{align}
		\begin{proof}
			We rewrite Equation \ref{eqn:matrixpredec} by re-indexing the last indices of both sequences. We now have that 
			\begin{align}
				\sum_{j = 0}^{n} s(n,j) F(j) = n!.
			\end{align}
			The infinite lower triangular matrix of $s(n,k)$ multiplied on vector $\vec{F}$ is exactly the sum derived. The equations
			\begin{align}
				\hat{s}\vec{F} &= \vec{f} \nonumber
			\end{align}
				and
			\begin{align}
				\vec{F} &= \hat{S} \vec{f} \nonumber
			\end{align}
		immediately follow, using the inverse of $\hat{s}$ being $\hat{S}$ \eqref{eqn:stirlinginvs}.
		\end{proof}
	\end{corollary}
	Corollaries \ref{cor:auseful} and \ref{cor:lint} show that Equation \ref{eqn:rfubfinal} provides intermediate transformations between strong and weak ordering. The intermediate sequences have combinatorial interpretation. Matrices $\hat{S}$ and $\hat{s}$ are lower triangular. Therefore, the relations of Corollary \ref{cor:lint} hold for finite cases, specifically for truncations of $\vec{f}$, $\vec{F}$, $\hat{s}$, and $\hat{S}$. 
	\section{Modular periodicity}
 Below, we determine eventual modular periodicity for $F(n)$ and $F_{r}(n)$ via the perspective of exponentially generated sequences and transformations of such sequences that preserve their structure. In addition, we set an upper bound for the period. The bound is the Carmichael function of the modulus. The case of odd modulus $K$ has a period $\lambda(K)$, which we prove directly by Asgari and Jahangiri's period formula \cite{cc:asgari}.
	\subsection{Fubini numbers modulo $K$}
	\begin{corollary}\label{thm:fubinimodk}
		The Fubini numbers are eventually periodic modulo $K = p_{1}^{R_{1}} p_{2}^{R_{2}} \cdots p_{N}^{R_{N}}$, with maximum possible period $\lambda(K)$. Periodicity holds for $n \geq \max(R_{1}, R_{2}, \ldots {R_{N}})$.
		\begin{proof}
			Consider the first relation of Corollary \ref{cor:lint}, which is
			\begin{align}
				\vec{F} = \hat{S}\vec{f}\nonumber
			\end{align}
			The entries of $\vec{f}$ are $n!$. Clearly, $n! \pmod{K}$ is zero for $n \geq K$. The relation stands in a simplified form modulo ${K}$ as
			\begin{align}
				F(n) \equiv \sum_{k = 0}^{K - 1} S(n,k) k! \pmod {K}\label{eqn:stirlingsecondmod}.
			\end{align}
			Fubini numbers modulo $K$ are written as a finite sum of $S(n,k)$ with coefficients independent of $n$. Each $S(n,k) \pmod {K}$ term contributes sums of exponential dependence in $n$ via Lemma \ref{lem:sskexplicit}. The weighted sum of modular exponentiations $j^{n}$, with $j\in \mathbb{N}_{0}$, is eventually modular periodic. Let $R = \max(R_{1}, R_{2}, \ldots {R_{N}})$ given $K = p_{1}^{R_{1}} p_{2}^{R_{2}} \cdots p_{N}^{R_{N}}$. The Carmichael function bounds the period of modular exponentiation since
			\begin{align}
				j^{R + \lambda(K)} \equiv j^{R} \pmod {K}.
			\end{align}
			All the exponentials $j^{n}$, each with fixed coefficients, must enter periodicity by $n = R$. Therefore, their sum is modular periodic for $n 
			\geq R$. The longest possible period of $F(n) \pmod {K}$ is $\lambda(K)$ because $\lambda(K)$ is the $\lcm$ of all periods of exponentiation modulo $K$. 
		\end{proof}
	\end{corollary}
	\subsection{Exact Carmichael periodicity for odd $K$}
	\begin{corollary}
		For $F(n) \pmod {K}$, if $K$ is odd, then the eventual modular period of the sequence is $\lambda(K)$.
		\begin{proof}
			According to Asgari and Jahangiri \cite[Thm.\ 10]{cc:asgari}, the period for this case is $\lcm(\varphi(p_{1}^{R_{1}}),\\ \varphi(p_{2}^{R_{2}}), \ldots, \varphi(p_{N}^{R_{N}}))$ where $p_{i}^{R_{i}}$ are the prime power factors of $K$. When $m$ is an odd prime power, we have $\varphi(m)=\lambda(m)$. The result is the following expression for the eventual period $q$:
			\begin{align}
				q &= \lcm(\lambda(p_{1}^{R_{1}}), \lambda(p_{2}^{R_{2}}), \ldots, \lambda(p_{N}^{R_{N}})).
			\end{align}
			Finally, we apply the recurrence for $\lambda(n)$ \eqref{eqn:carmichaelrecur} to find that
			\begin{align}
			q&= \lambda(K).
			\end{align}
		\end{proof}
	\end{corollary}
	\subsection{Extension to $r$-Fubini numbers}
	\begin{corollary}
		The $r$-Fubini numbers are eventually periodic modulo $K$ for fixed $r$, with maximum period $\lambda(K)$. The sequence's periodicity holds when $n \geq r - 1 + R$. Here $R = \max(R_{1}, R_{2}, \ldots {R_{N}})$ given $K = p_{1}^{R_{1}} p_{2}^{R_{2}} \cdots p_{N}^{R_{N}}$.
		\begin{proof}
			The Fubini numbers $F(n)$ are eventually modular periodic by Theorem \ref{thm:fubinimodk}. Next, observe the operation $r!\binom{E}{r} E^{-r}$ on $F(n)$ to generate $F_{r}(n) = r!\binom{E}{r} E^{-r} F(n)$, per Equation \ref{eqn:binomform}. The operation preserves the structure of $F(n)$ as a weighted sum of exponentials of $n$. The sum modulo $K$ has a maximum period of $\lambda(K)$. Shift operators might delay the periodic onset for some terms, so we add $r - 1$ to the onset index. Consider that $F_r(n)$ includes terms with argument shifts of $F(n)$. Some $F(n - a)$ term in $F_{r}(n)$ only has argument $R$ (ending its initial aperiodic procession) when $n = a + R$. The largest such $a$ is $r - 1$ according to Equation \ref{eqn:rfubfinal}. Note $s(b,0) = 0$ when $b > 0$.
		\end{proof}
	\end{corollary}
	\section{Remarks}
	\subsection{Analogy between ordered and unordered Bell numbers}
	\begin{remark}
		Consider the formula for $r$-Bell numbers given by Nyul and Keresk\'enyi-Balogh \cite[Thm.\ 4.2]{cc:nyul}, rewritten with $n,r$ in our notation it reads
		\begin{equation*} B_{r}(n) = 		
		\sum_{j = 0}^{r} s(r,r - j) B(n - j).
		\end{equation*} The formula for $r$-Fubini numbers given in this work (\ref{eqn:rfubfinal}) is remarkably similar. The proof given here should apply to $r$-Bell numbers with some adjustment. Explaining this similarity could yield insight into other problems susceptible to operator counting.
	\end{remark}
	\section{Acknowledgments}
	William Gasarch gave great advice for generalizing from an inspired example to $H_{r}(n)$ and supported writing. The editor of the journal and referees contributed valuable feedback on writing. Alex Leonardi, Deven Bowman, Kevin Flanary, Nathan Constantides, and Elan Fisher gave crucial review or aided in brainstorming.
	
	\bigskip
	\hrule
	\bigskip
	\noindent 2020 {\it Mathematics Subject Classification}: Primary 06A05; Secondary 11B50.
	
	\noindent \emph{Keywords:} Carmichael function, constrained weak ordering, exponential sum, Fubini number, modular periodicity, ordered Bell number, $r$-Fubini number, $r$-horse number, shift operator, signed Stirling number of the first kind, Stirling number of the second kind, weak ordering.

	\bigskip
	\hrule
	\bigskip
	
	\noindent (Concerned with sequence
	\seqnum{A000670}, \seqnum{A002322}, \seqnum{A008275}, \seqnum{A008277}, \seqnum{A232473}, 
	and \seqnum{A232474}.)
	
	\bigskip
	
	\bigskip

	\noindent

	\bigskip
	
	\bigskip
	
	\noindent
	
	\vskip .1in

\end{document}